# Parameter Estimation of Noise Corrupted Sinusoids
## An Engineering Method

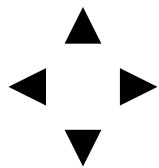

Francis J. O'Brien, Jr., Ph.D. and Nathan Johnnie, P.E.
Naval Undersea Warfare Center, Division Newport
Undersea Warfare Combat Systems Department
Code 2531, Bldg. 1259
1176 Howell St.
Newport, RI 02841-1708

**May 2, 2011**

## Abstract

Existing algorithms for fitting the parameters of a sinusoid to noisy discrete time observations are not always successful due to initial value sensitivity and other issues. This paper demonstrates the techniques of FIR filtering, Fast Fourier Transform, circular autocorrelation function, and nonlinear least squares minimization as useful in the parameter estimation of amplitude, frequency and phase, exemplified for a low-frequency time-delayed sinusoid describing simple harmonic motion. Alternative means are described for estimating frequency and phase angle. An autocorrelation function for harmonic motion is also derived.

## 1. General Purpose

The sine wave or sinusoid is a well known mathematical function. Periodic functions with real or imaginary roots are encountered in many fields to describe smooth repetitive motion in time and space. Oscillation is modeled by linear second order homogeneous differential equations with constant coefficients derived from the equation, $\frac{d^2x}{dt^2} + \omega^2 x = 0,$ where $t$ represents time and $\omega = 2\pi f$ is the angular frequency in radians per second such that $\frac{1}{f} = \frac{2\pi}{\omega}$ is one period. The general solution is simple harmonic motion (undamped oscillation), $x(t) = A\sin(\omega t + \varphi)$, for noise free systems where $\varphi$ is the initial phase angle. Finney and Thomas (Chapter 16) provide a brief but complete derivation of $x(t)$. More detailed definitions are given below in Notation.

Since sinusoids are basically graphs of wave forms, then any phenomena having periodic behavior or wave characteristics can be represented by trigonometric functions or (approximately) modeled by sinusoids. A large number of fields employ sinusoids including physics, electrical engineering, mechanical engineering, signal processing, economics/finance, neural networks, and others. This includes many simple actions such as a weighted spring, pendulum, motion of an engine's piston-crankshaft, ocean tides, blood pressure in the heart, sounds, air temperature, and many other phenomena either created by nature or man-made machines/systems.

Modeling simple harmonic motion by parameter estimation of sinusoids in the presence of noise is the focus of the present effort. This paper will provide a useful method to obtain parameter estimates of the three cardinal terms of a simple sinusoidal signal, $x(t) = A\sin(\omega t + \varphi)$, and the autocorrelation function (ACF) of the sinusoid for a noise–corrupted finite time series data set. An example for a time-delayed low frequency distribution of synthetic data is given.

Existing methods such as the IEEE Standard 1057 (IEEE–STD–1057) provide algorithms for fitting the parameters of a sine wave to noisy discrete time observations. The engineering literature indicates that the standard nonlinear curve fitting approaches are not always successful due to initial value sensitivity and other issues. For these reasons the authors have derived the present method for a range of frequencies, phases and sample sizes that is not dependent on initial value specification.

A discrete nonparametric circular autocorrelation function (ACF) (Wald and Wolfowitz, 1943) is used to initially describe the noisy data $X(t)$. This measure provides a normalized ACF for periodic data with a well defined symmetric cosine wave shape unlike the die–away appearance of a standard ACF for discrete data. A visual indication of periodic structure in the data can be determined by inspecting the pattern of the ACF values, although tests for periodicity based on a spectral analysis of a time series are normally used; e.g., The Fisher Test for the detection of hidden periodicities.

The sinusoidal model of $x(t)$ is characterized by a newly derived ACF calculated from the integral of a sine function. This ACF model provides two possible solutions to characterize the sinusoid of simple harmonic motion. Each was analyzed and the best one selected. Appendix I shows the logic.

The initial noisy time series data set $X(t)$ must be reduced to a manageable mathematical function in order to model the data parametrically. To accomplish this, the time series must first be smoothed by means of a moving average (MA) filter or other filtering techniques, well known to those skilled in the art. These initial operations are needed to create an approximating function that attempts to capture important periodic structures in the data, while leaving out noise or other fine–scale structures/rapid phenomena of the erratic ups and downs of the series.

To smooth the data the authors employ a statistical filter known as a moving average (MA-$k$) which is a type of discrete-time finite impulse response (FIR) filter. It is commonly used with time series data to smooth out short–term fluctuations and highlight longer–term trends or cycles/periodicities and is also similar to the low–pass filter. MA filters of size MA-5 and MA-10 are usually sufficient to model noisy data for short-term and long-term fluctuations.

**NOTE:** The discrete FIR filter was judged superior to a standard DSP smoothing filter in providing a smooth periodic curve.

As a specific example of the method process, the following EXCEL charts exemplify a 5–and 10–point moving average (black line) for 100 noisy data points (blue line) derived from a periodic function with amplitude of 2, frequency of 0.05 Hz, phase of 35 degrees $\left( \varphi = \frac{7}{36}\pi = 0.6109 \, \text{rads.} \right)$. The noise level was set in MATLAB at 0.5*randn (1/2 of the original signal is Gaussian noise). See Graph A below.

The 5-point MA shows the best goodness of fit in this case. The authors are unable to prescribe the optimal MA smoothing parameter, as this is an unsolved issue in data processing[1]. Root mean square (RMS) error analysis is the most common method to judge effectiveness of smoothing.

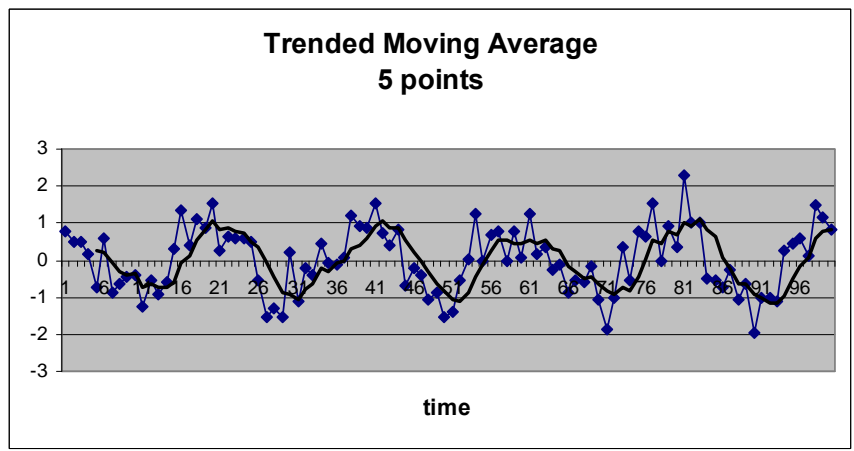

---

[1] As a general guideline in the application of moving averages for detecting periodicities, the MA parameter $k$ should be selected to be an even multiple of the period $T$. This suggests that frequency should be obtained first by the Fast Fourier Transform (FFT) or reading the autocorrelation functions described below as a cosine function graph; this is illustrated below in Figures 1 & 3.

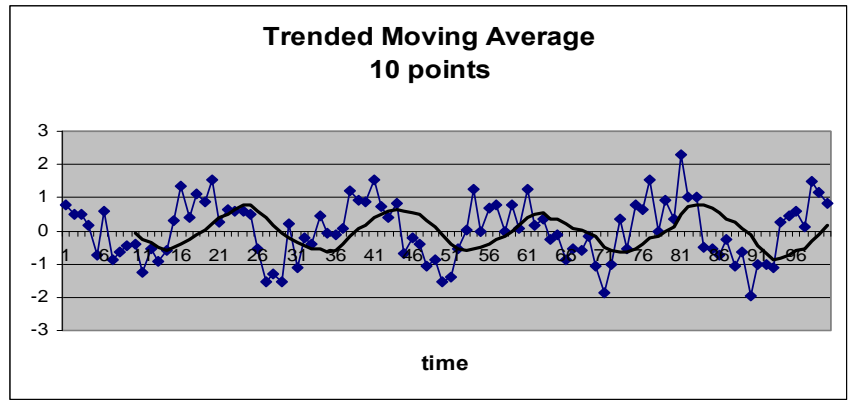

Graph A. 5– and 10–point moving averages (MA-*k*) for noisy data.

The performance analysis of MA filtering showed several key factors germane to the engineering method:

- A raw data mean of about 0.01 indicates a stationary process
- The structure of the MA graph indicates a periodic signal which can be further corroborated by Fisher's test or alternative means
- The original data also showed amplitude close to 2.1 (MA-5)
- A Fast Fourier Transform (FFT) on the raw data will show the fundamental frequency to be 0.05 Hz. See next graph (Graph B) for this level of noise; higher noise levels may not give a clean reading of *f*. This provides the period, $T = \frac{2\pi}{\omega} = 20$. However, either moving average filter provides a reading of about $T = 20$. The autocorrelation functions (Figures 1,3 below) confirm this finding.
- Since phase cannot be calculated directly from the FIR filter a nonlinear least squares minimization search algorithm is used to calculate $\varphi$ as described below (Eq. 8). Alternatively, we can estimate $2\pi$ (at second crossover of MA-10) at the 18$^{th}$ position, $t_{18}$. Since the phase displacement in absolute magnitude is $\Delta t = \frac{1}{f} - t_{2\pi} = \frac{\varphi}{360 f}$ (degrees), we obtain an initial estimate of phase as $\Delta t = 2$ or $\varphi = 36° = \frac{\pi}{5} = 0.6283$ rad. This estimate is derived from the relation, $\omega t + \varphi = 2\pi$ or, $t_{2\pi} = \frac{2\pi - \varphi}{\omega} = T - \Delta t$. Because the phase is positive then $\Delta t$ is negative (time delay) from the relation, $\text{Arcsin}[A\sin(\omega t + \varphi)] = 0$, or $t_0 = \frac{-\varphi}{\omega} \Rightarrow \Delta t = -2.0$. The total error of estimate for $\varphi$ is about 2.9%

Thus, crude but fairly accurate estimates are possible with the MA filters for this level of noise, and the results provide a preliminary solution for the parameters *A*, *f* & $\varphi$ of the noisy sinusoid, $x(t) \approx 2.1\sin\left(\frac{2\pi}{20}t + \frac{\pi}{5}\right)$.

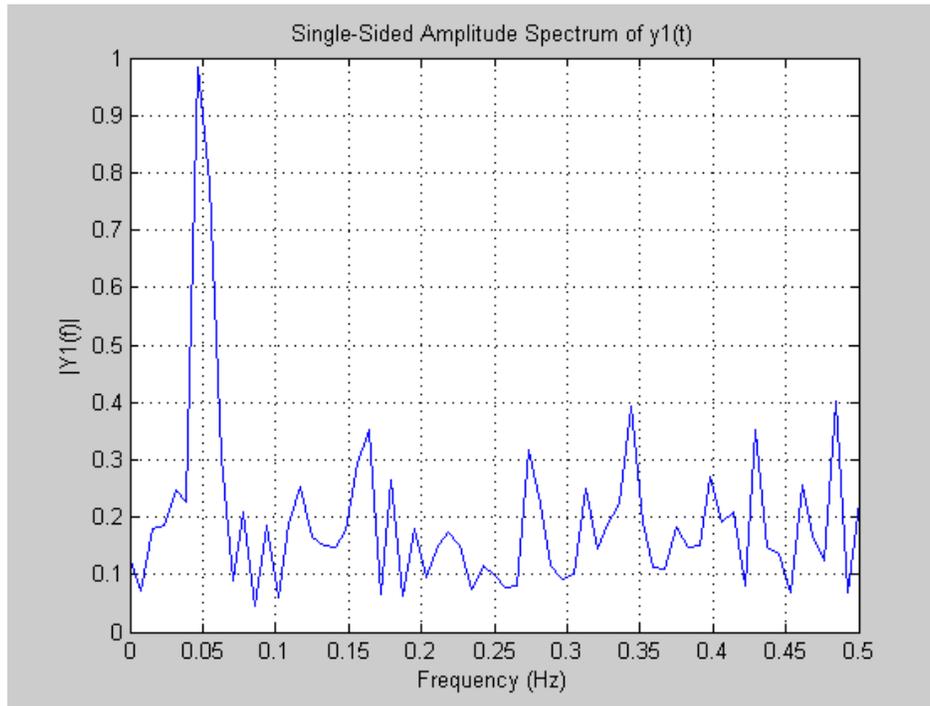

Graph B. Fast Fourier Transform (FFT) on data with 50% Gaussian noise showing fundamental frequency $f = 0.05$ Hz.

An initial assessment of the noisy data set $X(t)$ is made to determine if any useful signal information exists in the data or they are "pure noise" (random). To this end the non–parametric Wald–Wolfowitz Runs Test method described in US Patent 7, 103, 502, Sep. 5, 2006, awarded to F. J. O'Brien, Jr., is implemented to provide a rapid decision of "signal" or "noise". A decision of "noise" indicates no useful signal information can be extracted at which point processing ceases for that data set.

A second screening test with the same objective is performed with the modified Wald—Wolfowitz discrete circular ACF. We have derived an alternative procedure for interpreting the Wald–Wolfowitz circular correlation outcomes in the presence of noise which improves the accuracy of the method[2]. This discrete ACF was shown to be symmetric on $t \geq 0,$ a unique characteristic of a serial correlation function. The implication of the fold–over symmetry property is that a complete ACF can be created in only ½ the steps ordinarily needed to perform an ACF for the maximum number of lags. In addition, this nonparametric ACF is applicable for data/noise that is either Gaussian or non–Gaussian.

The purpose of the 2–gate screening is to discard inputs which are essentially pure noise. The present work assumes that the data set $X(t)$ contains some degree of "signal" in the data structure in order to perform useful parametric modeling of the sinusoid of simple harmonic motion, $x(t) = A\sin(\omega t + \varphi)$.

---

[2] The procedure is a Navy Invention Disclosure, and is subject to non disclosure at the present time.

## 2. Background

**NOTATION**

- $t$ is time, discrete or continuous
- $N$ is sample size
- $\tau$ is the time lag
- $x(t)$ is hypothetical sinusoidal model, $A\sin(\omega t + \varphi)$ with no additive vertical shift factor
- $X(t)$ is noise corrupted input time series
- $\eta$ is is noise level of data set $X(t)$, Gaussian or non-Gaussian
- $A$ = amplitude or maximum of the function relative to the horizontal time axis
- $\omega = 2\pi f$ is the angular frequency which specifies how many oscillations occur in a unit time interval, in radians per second (in Hertz). Also called radian frequency (rad/sec)
- $T = \dfrac{2\pi}{\omega} = \dfrac{1}{f}$ is time duration of one period or the time required for one complete cycle, revolution or oscillation
- $f = \dfrac{\omega}{2\pi}$ is the fundamental frequency or the number of complete periods per unit time in Hz
- $\varphi$ = initial phase angle of the motion, or phase offset, expressed in radians & defined as the fractional part of the period which has elapsed at any instant since the object passed through its central position in the positive direction, $-\pi \leq \varphi < \pi$
- $\Delta t = \dfrac{\varphi}{360 f}$ (degrees) $= \dfrac{\varphi}{\omega}$ (radians) is the time delay (or time ahead) factor. Setting $t \to t \pm \Delta t$ in the sinusoid will create the primitive function, $A\sin(\omega t)$. $|\Delta t|$ is also called the "phase displacement".

### General Integral Solution of Autocorrelation Function

Let the autocorrelation function (ACF) for continuous measurements in finite time $t$ for a time lag $\tau$ be defined as:

$$R_{xx}(\tau) = \frac{1}{T}\int_0^T x(t)x(t+\tau)dt, \quad (t,\tau) \geq 0 \qquad (1)$$

where, $T$ is defined above in Notation. An important property of an ACF is that periodic data have periodic ACFs with the same frequency.

This function must be normalized appropriately to provide $\pm 1$ values. We consider different solutions to Eq. (1) for characterizing sinusoidal functions of the following form:

$$\begin{cases} x(t) = A\sin(\omega t + \varphi) \\ x(t+\tau) = A\sin[(t+\tau)\omega + \varphi] \end{cases} \tag{1a}$$

A general solution to the ACF of Eq. (1) for sine functions is achieved by evaluating the elementary indefinite integral,

$$\int \sin(ax+b)\sin(cx+d)dx \tag{2}$$

Signal processing solutions of (2) are provided when $a = c$; i.e., sinusoids with the same value of $\omega t$. Two models can be created by (2) for (1a) functions depending upon the term of integration, $dt$ or $d\varphi$.

Eq. (2) is evaluated first by expansion in accordance with the identity $\sin A \sin B$ and a change of variable which gives the result,

$$\int \sin(ax+b)\sin(ax+d)dx = \frac{x}{2}\cos(b-d) - \frac{\sin(2ax+b+d)}{4a} \tag{3}$$

**NOTE:** if the sine term is replaced by cosine in the integrand, it differs from (3) only by being a sum of the two terms rather than a difference. Other indefinite forms are found in Gradshteyn & Rhyzik, Sect. 2.542.

Eq. (3) provides the means to derive a general definite integral solution for arbitrary limits $(u,v)$ by identity $\sin A - \sin B$:

$$\int_u^v \sin(ax+b)\sin(ax+d)dx = \frac{v-u}{2}\cos(b-d) - \frac{\sin[a(v-u)]\cos[a(u+v)+b+d]}{2a} \tag{4}$$

This integral allows examination of the ACF for defined regions.

### 3. Description and Operation

Specific ACF Solutions

#### A. Full Model

Substitution for the limits $(u,v)$ provides specific signal processing solutions to Eq. (4) with $x = t, a = \omega, b = \varphi, d = \tau\omega + \varphi$.

- If $[u = 0, v = T]$, one ACF is:

$$R_{xx}(\tau) = \frac{1}{T}\int_0^T A\sin(\omega t + \phi) A\sin(\omega t + \omega\tau + \phi)dt$$

$$R_{xx}(\tau) = \frac{A^2}{2T}\left\{T\cos(\omega\tau) - \frac{\sin[\omega T]\cos[(T+\tau)\omega + 2\varphi]}{\omega}\right\} \tag{5}$$

$$0 \leq \tau \leq N-1.$$

**NOTE**: The authors have used this empirically normalized ACF with useful results. The following figure (Fig. 1) shows one typical run on a noisy sinusoid with Gaussian noise set at 0.5*randn. For this plot the inputs are $A = 1, f = 1/20\,\text{Hz}, \phi = 35°\,(0.6109\,\text{rads})$. The red line is the empirical sinusoid for the time series of 100 observations. The cyan line is the model (5) ACF. Note that the second maximum value (20) of either ACF gives the period $T$.

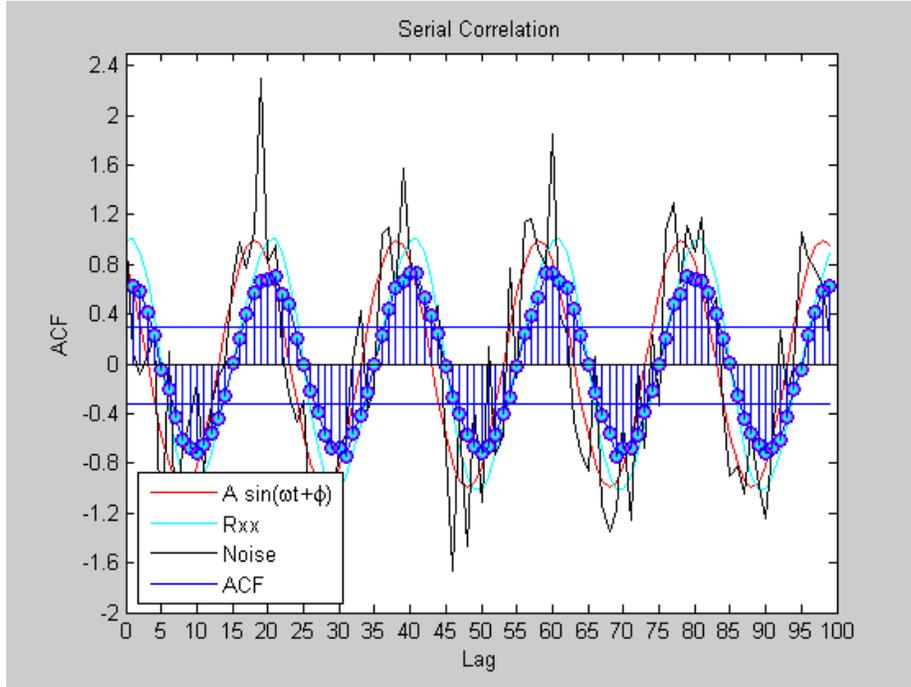

Figure 1. Normalized ACF, Eq. (5) (cyan colored line, labeled "Rxx") applied to noisy input function, $y = \sin\left(\frac{2\pi}{20}t + 35°\right) + \eta$ (black line, labeled "Noise"). Gaussian Noise $\eta$ is 0.5*randn. Sample size is $N$ = 100. False alarm rate is 0.001. The method judged this data set "periodic signal". NOTE: The 0–lag of circular "ACF" is omitted to declutter the 4–plot graph.

- The analytic normalization of (5) is provided next. Thus, as a second solution, we considered $[u = 0, v = 2\pi]$, for one period, to obtain the ACF:

$$R_{xx}(\tau) = \frac{1}{2\pi}\int_0^{2\pi} A\sin(\omega t + \phi)\,A\sin(\omega t + \omega\tau + \phi)\,dt$$

$$R_{xx}(\tau) = \frac{A^2}{2\pi}\left\{\pi\cos(\omega\tau) - \frac{\sin[2\pi\omega]\cos[(2\pi + \tau)\omega + 2\varphi]}{2\omega}\right\}$$
(5a)

To normalize this formula we use the definition $R_{xx}(\tau) = 1$ for $\tau = 0$, so that the normalizing constant $C$ is:

$$C = \frac{2}{A^2}\left[1 - \frac{\sin(2\pi\omega)\cos(2\pi\omega + 2\varphi)}{2\pi\omega}\right]^{-1}$$

to give a normalized autocorrelation function,

$$R_{xx}(\tau) = C\left\{\cos(\omega\tau) - \frac{\sin[2\pi\omega]\cos[(2\pi+\tau)\omega+2\varphi]}{2\pi\omega}\right\} \quad (5b)$$

where,

$$C = \frac{2}{A^2}\left[1 - \frac{\sin(2\pi\omega)\cos(2\pi\omega+2\varphi)}{2\pi\omega}\right]^{-1}.$$

Thus, for the data $A = 2, f = 1/20 \text{ Hz.}, \phi = 35°\ (0.6109 \text{ rads})$, we calculate the normalizing constant as $C \approx 1.465315522\left(\frac{2}{A^2}\right)$

If we ignore altogether the scaling amplitude factor, $\frac{A^2}{2}$, then (5b) can be stated in a form amenable to machine calculation for sample size $N$ with $\tau \geq 0$:

$$R_{xx}(\tau) = \frac{1}{1 - \frac{\sin[2\pi\omega]\cos[2\pi\omega+2\varphi]}{2\pi\omega}}\left\{\cos(\omega\tau) - \frac{\sin[2\pi\omega]\cos[(2\pi+\tau)\omega+2\varphi]}{2\pi\omega}\right\}$$

(5c)

$-1 \leq R_{xx}(\tau) \leq 1,\ 0 \leq \tau \leq N-1.$

Eq. (5c) is the optimal modeling ACF of sinusoids of the form $A\sin(\omega t + \varphi)$ derived from noisy data.

The normalized ACF of Eq. (5c) was validated analytically by calculating values for $\tau$ at known points of the max/min peaks and troughs of the noise free cosine function at $\left[0, \frac{\pi}{2}, \pi, 2\pi\right]$ which corresponds to $\tau = [0, 5, 10, 20]$ for $f = 0.05$ Hz. on one period. The following table shows that the actual values from Eq. (5c) are as predicted.

| $\tau$ | $R_{xx}(\tau)$ Eq. (5c) |
|---|---|
| 0 | +1.00 |
| 5 | −0.017 |
| 10 | −1.00 |
| 20 | +1.00 |

**NOTE**: The formulation of Eq. (5c) treats all terms of the sinusoid $(A, \omega, \varphi)$ as unknown random variables. Eq. (5c) has been tested empirically on noisy data to determine its properties to describe noisy data with hidden periodicities. However, Eq. (5c) cannot be used to estimate the frequency in closed form from a noise–corrupted time series data set $X(t)$. Consequently, a second model was considered for calculation of $f$.

### B. Reduced Model

- If the ACF of (1) is integrated in terms of a random phase $\varphi$ with assumed constant amplitude and $f$, a simpler and perhaps more practical modeling solution results. This calculation, well known to those skilled in the art, results in:

$$R_{xx}(\tau) = \frac{1}{2\pi} \int_0^{2\pi} x(t)x(t+\tau)d\varphi$$

$$R_{xx}(\tau) = \frac{1}{2\pi} \int_0^{2\pi} A\sin(\omega t + \varphi)A\sin(\omega t + \omega\tau + \varphi)d\varphi \qquad (6)$$

$$R_{xx}(\tau) = \frac{A^2}{2}\cos(\omega\tau)$$

**NOTE:** Eq. (6) gives the ACF of a sinusoid, $A\sin(\omega t + \varphi)$, with the same frequency which can be estimated at a time $t$. However an FFT of noisy data provides $f$ readily, but (6) can be used a verification check or for models that consider phase as a random variable and frequency/amplitude as fixed.

**NOTE**: Analysis shows that Eq. (6) also results if cosine terms replace sine terms in the integrand.

**NOTE**: Eq. (6) is normalized by ignoring the scaling factor, $\frac{A^2}{2}$, since by definition, the cosine function has the property, $-1 \leq \cos(\omega\tau) \leq +1$.

**NOTE:** Although different theoretical models, (5c) and (6) are highly correlated numerically, differing on average about 2% in value at any lag. Fig. 3 below exhibits the similarity in structure of each ACF model.

The ACF of Eq. (5c) is preferable since it models all parameters of the sinusoid. Eq. (6) provides a symmetric cosine even function while (5c) is not strictly symmetric.

The following table shows the values of each ACF for 1/2 period for $T = 20, f = 0.05$ Hz.:

| Tau $\tau$ | ACF Eq. 5c | ACF Eq. 6 |
|---|---|---|
| 0 | 1 | 1 |
| 1 | 0.9457 | 0.9511 |
| 2 | 0.7989 | 0.8090 |
| 3 | 0.5739 | 0.5878 |
| 4 | 0.2922 | 0.3090 |
| 5 | −0.0172 | 0 |
| 6 | −0.3254 | −0.3090 |
| 7 | −0.6017 | −0.5878 |
| 8 | −0.8191 | −0.8090 |
| 9 | −0.9564 | −0.9511 |
| 10 | −1 | −1 |

The simple linear correlation for these values is about 1.00.

**NUMERICAL EXAMPLES**

These examples assume knowledge of values of $A, f, \varphi$ which are estimated as shown in the next section.

- Figure 2 indicates that the theoretical normalized cosine ACF of Eq. (6) — $\cos(\omega\tau)$ — models exactly as predicted from the calculation. Eq. (5c) performs at about the same level. The theoretical ACF plot (in cyan) is not shown as it lays beneath the discrete circular ACF values (indicated by blue dots). The "Noise" curve is the input sine function since the sinusoid is noise free. The input sine function was $y = 2\sin\left(\frac{2\pi}{20}t + 35°\right)$ with an amplitude of $A = 2$, $f$ = frequency of 1/20 Hz and phase $\varphi$ of 35° (0.6109 radians) and time delay, $\Delta t$ of $\left|\frac{\varphi}{\omega}\right| \approx 1.94$. The discrete circular ACF was computed for this circular function as an initial assessment of periodicity and as part of the screening process. The false alarm rate (FAR) is 1%.

- Figure 3 shows Eq. (5c) and Eq. (6) applied to a noisy sinusoid of the form, $y = \sin\left(\frac{2\pi}{20}t + 35°\right) + \eta$. Gaussian noise (0.5*randn) $\eta$ was added to $y$. The ACFs (5c) and (6) appear indistinguishable (magenta and green plots)  Error bounds on the discrete ACF (in blue) are indicated by the solid blue lines parallel to the lag–axis.
- Figure 4 shows an example of "pure noise" rejected by the two-stage screening process. The FAR was run at a higher level (0.001). It also shows why it must be rejected. The degree of error (difference between ACF (blue dots) and the modeled ACF (cyan line)) is so large that no meaningful estimation of parameters is possible.

- Many other calculations and tests were performed at higher noise levels. The data sets were first screened with the rapid Runs Test of Pat. '502. Distributions that passed this test (deemed "signal") were processed by the present method. The next gate consisted of testing the circular discrete ACF to determine if the distribution was signal or noise.  The model tracked the noisy periodic signal but with less goodness of fit as would be expected for highly noise dominated data. Lower FARs (e.g., 0.001 or less), when used in conjunction with the screening procedure,  can improve efficiency.

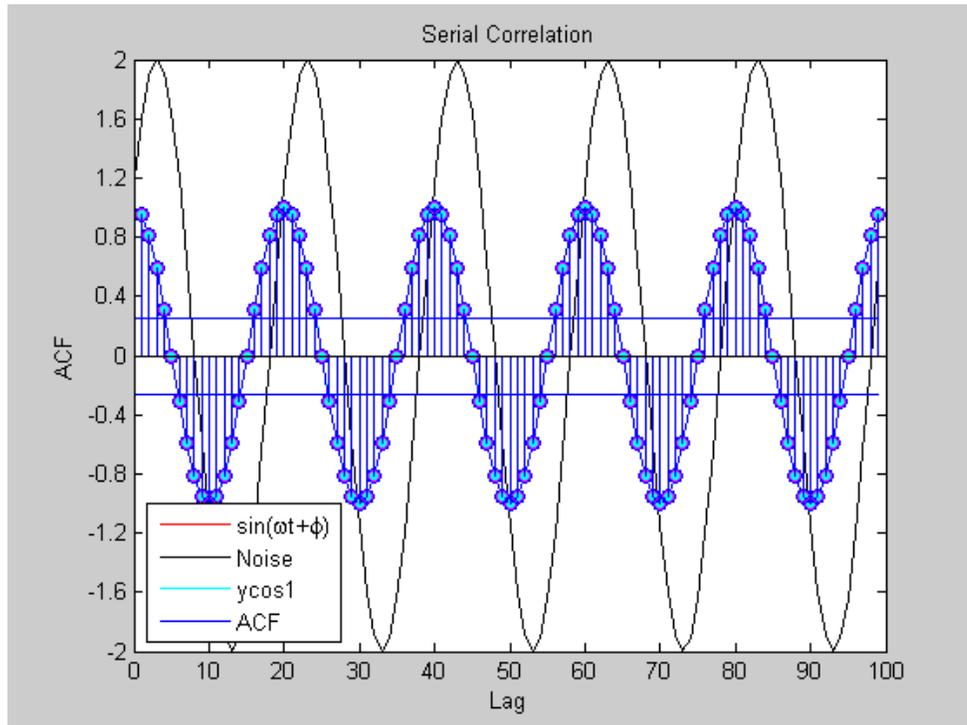

Figure 2. Theoretical cosine ACF $\cos(\omega\tau)$ (cyan colored line, labeled "ycos1", applied to noise–free input function, $y = 2\sin\left(\dfrac{2\pi}{20}t + 35°\right)$ (black line, labeled "Noise"). Sample size is $N = 100$. False alarm rate is 0.01. The method judged this data set "periodic signal"  The symmetry around the 50$^{th}$ lag $(N/2)$ is characteristic of this ACF for discrete data.

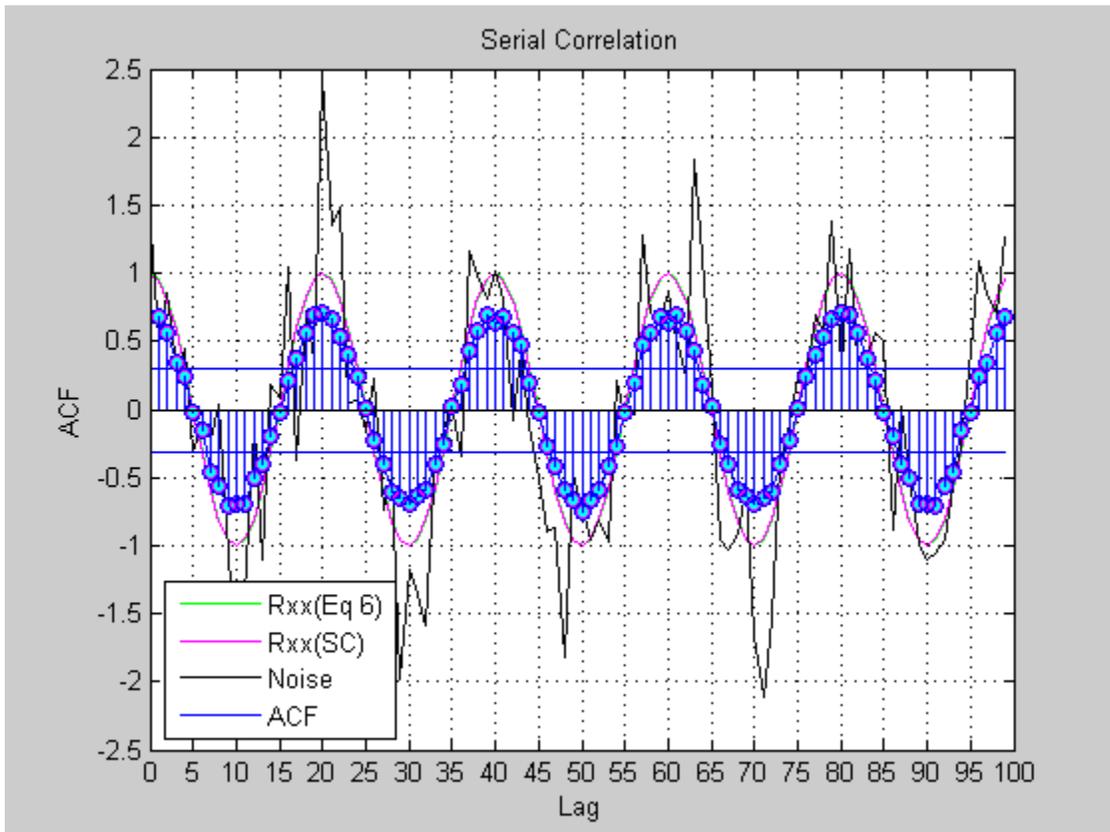

Fig. 3. Theoretical cosine ACF Eq. (5c) (magenta colored line, labeled "$R_{xx}$"). This ACF appears to be the same as Eq. (6) (green colored line) applied to noisy input sine function $y = \sin\left(\dfrac{2\pi}{20} t + 35°\right) + \eta$ (black line, labeled "Noise"). Noise is 0.5*randn. Sample size is $N = 100$. The method detected "periodic signal" for these data. False alarm rate is 0.01. NOTE: The 0–lag of "ACF" is omitted to declutter the 4–plot graph. The line labeled "ACF" is the circular correlation function for the sinusoid.

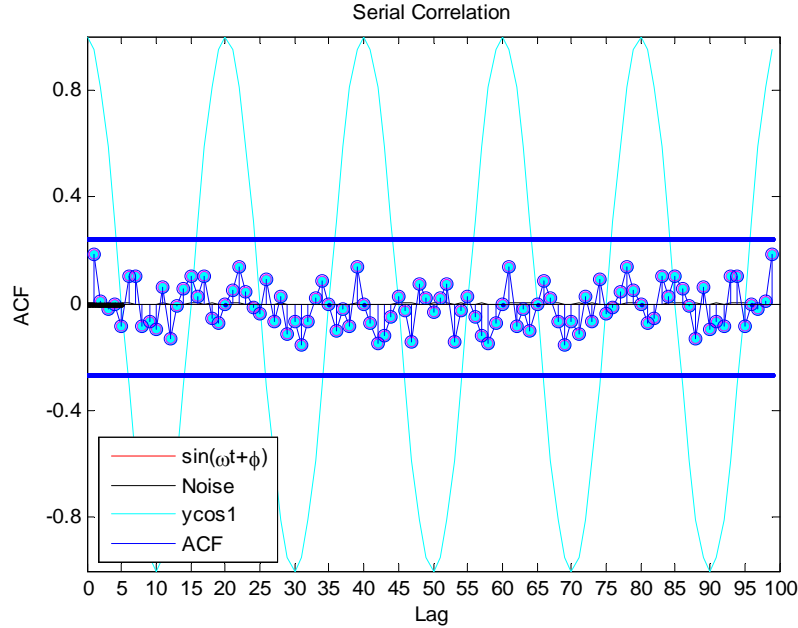

Figure 4. Example of "pure noise" rejected by screening process. The raw data and sine function are not plotted. Noise is 80*randn. Sample size is $N = 100$. The FAR was 0.001. Obviously no modeling is possible for these random data.

## Estimation of Parameters $A, f, \varphi$

The FIR smoothing filter and the nonlinear estimation algorithms will characterize a simple sinusoidal function, $A\sin(\omega t + \varphi)$, from a time series of noisy data. This is the main objective of the present engineering method.

### Estimating amplitude $A$

Known range-based methods exist for estimating the amplitude $A$, so we assume a fairly accurate estimate of $A$ exists from the noisy data set. However, the method will suffice based on the smoothing FIR filter selected such as described earlier.

### Estimating frequency $f$

The numerical value of the frequency can be derived in closed form from Eq. (6) using the arccosine function, $\cos^{-1}(x), 0 \leq x \leq 1$. For any value of the ACF of a noisy sinusoid we can compute by (6):

$$R_{xx}(\tau) = \cos(\omega\tau) = \cos(2\pi f\tau).$$

Then,
$$\cos^{-1}[R_{xx}(\tau)] = \cos^{-1}[\cos(2\pi f\tau)] = 2\pi f\tau$$

$$f = \frac{\cos^{-1}[\cos(\omega\tau)]}{2\pi\tau}, \tau \neq 0. \quad (7)$$

Example: To illustrate, assuming Fig. 2 represents the best fitting sinusoid of the data $(A = 2, f = .05, \varphi = 0.6109)$, the tau value $\tau$ at 2 is $R_{xx}(\tau = 2) = \cos(2\omega) = 0.8090$ so that

$$f = \frac{\cos^{-1}(0.8090)}{2\pi\tau} = 0.05\,\text{Hz}.$$

which corresponds to the known $f$ for this hypothetical sinusoid Since $\frac{1}{f} = T$, this translates into one period, $T = 20\,\text{sec}$.

**NOTE:** If the value of $R_{xx}(2) = 0.7989$ from Eq. (5c) is used in (7) in place of 0.8090 from (6) as the best estimate, the results are almost identical, $T = 19.47$.

Typically, the FFT is the preferred method of obtaining $f$ from noisy data but (7) can be used as a confirmatory check on the calculation as shown above in Graph B. As noted the FFT is sensitive to the noise level. The maximum frequency is selected as the estimate of $f$ on the FFT plot when this is possible. Other methods can be used to obtain frequency $f$ as noted in Alternatives.

Estimating phase angle $\varphi$

Phase is more difficult to obtain accurately from noisy data. Sometimes only a crude estimate is available which must be refined in a nonlinear minimization least squares search as described below. The authors estimate phase angle $\varphi$ empirically. The object is to estimate $\varphi$ given knowledge of $A = 2$, $f = 0.05$ Hz and $\omega = 2\pi f = 0.3142$ as previously established, in the hypothetical modeling relation, $x(t) = 2\sin\left(\frac{2\pi}{20}t + \phi\right)$ which we assume arises as estimates of parameters of the sine wave derived from noisy data.

Phase angle measurement is performed preferably in the standard manner on a right triangle (in radians). The technique is described briefly in Finney and Thomas (chapter on differential equations and harmonic motion, pp. 1035ff.).

This technique is demonstrated for the noise free data in Fig. 2 as illustrated in the following diagram which is assumed for illustrative purposes to represent the derived sinusoid.

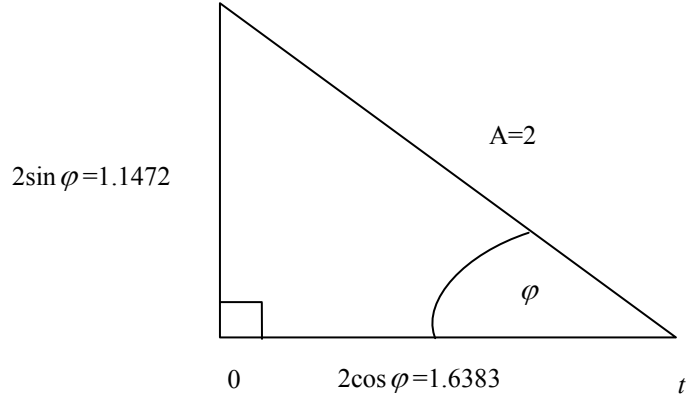

Phase angle $\phi$ measurement, Fig. 2.

For the vertical axis we obtain the quantity $x(t) = 2\sin\left(\frac{2\pi}{20}t + \phi\right) = 2\sin(\varphi) = 1.1472$ at the origin $t = 0$ of the sine wave. The other leg of the triangle is by subtraction using the Pythagorean Theorem, $\sqrt{4-(A\sin\varphi)^2}$, giving $2\cos(\varphi) = 1.6383$. Then,

$$\tan(\phi) = \frac{A\sin\varphi}{A\cos\varphi} = \frac{1.1472}{1.6383} = 0.7002$$

which leads to,

$$\varphi = \tan^{-1}\left(\frac{2\sin\varphi}{2\cos\varphi}\right) = \tan^{-1}(0.7002) = 0.6109,$$

the known input value of the phase angle. The $\pm$ sign of $A\sin\varphi$ determines if phase $\varphi$ lies above (positive) or below (neg.) the time axis on time $t \geq 0$ as well as whether the time delay, $\Delta t = \frac{\varphi}{\omega}$, is positive or negative.

For noise corrupted data the authors proceed in a manner based on the estimated parameters of the best fitting smoothed sine curve, $A\sin(\omega t + \varphi)$, to the observed data. In this procedure the object is to find a nonlinear least squares minimization solution to the following objective function on $t \geq 0$:

$$\sum_{t=0}^{T}[X(t) - A\sin(\omega t + \varphi)]^2 \rightarrow \min \qquad (8)$$
$$0 \leq t \leq T; -\pi \leq \varphi < \pi$$

$X(t)$ is the input signal + noise data. In (8) only $\varphi$ is assumed unknown since $A$ & $f$ have been obtained as explained previously.

The complete search in (8) ranges from 0 to $T$ and $\varphi$ varies from $-\pi$ to $\pi$ in increments of hundredths. As a second refinement, once phase is estimated to hundredths accuracy, a thousandths search accuracy can be achieved on the interval of the initial estimate of $\varphi$. The initial estimated value of $\varphi \approx 36°$ given earlier from the MA analysis could be used to speed up the algorithm.

Based on the FIR filter estimates, the following standard formulas in degrees and radians provide a quick calculation of $\varphi$ and acts as a check on the results of the nonlinear estimation of Eq. (8):

$$\left. \begin{array}{l} \varphi° = \dfrac{T - t_{2\pi}}{T}(360) = \dfrac{\Delta t}{T}(360) \\ \varphi_{rad} = \dfrac{2\pi(T - t_{2\pi})}{T} = \omega \Delta t \end{array} \right\} \quad (9)$$

If $\varphi = 0$, then $\dfrac{1}{f} = t_{2\pi}$, and Eq. (9) will return the correct value of 0 (no phase).

Eq. (9) can be generalized to the following formulations which allow phase to be computed anywhere on the sinusoid depending upon the measurements available,

$$\left. \begin{array}{l} \varphi_{rad} = \dfrac{T - t_{2\pi}}{t_{k\pi} - t_{2\pi}}(k\pi - 2\pi), \quad (k \neq 2) \\ \varphi° = \dfrac{\varphi_{rad}}{\pi} 180 \end{array} \right\} \quad (10)$$

which reduce to $\omega \Delta t$ upon simplification.

This formula follows from the definitions of $\varphi$ and time delay, $T - t_{2\pi} = \dfrac{\phi}{\omega}$ or $\varphi = \omega \Delta t$, for time-delayed sinusoids.

The result is the best fitting empirical sinusoidal curve to the noisy data. The periodic ACF of (5c) is then obtained for the solution, $A\sin(2\pi f t + \varphi)$.

Other methods for obtaining $\phi$ may be used, such as regression analysis based on the function (8). Another method, quite straightforward, can be used when a good empirical graph of $A\sin(2\pi f t + \phi)$ from noisy data is available such as the MA filters. Call that value $y$ so that $\varphi = \sin^{-1}\left(\dfrac{y}{A}\right) - \omega t$. For our data, this calculation shows, for selected $t = 1.8$, $\phi = \sin^{-1}\left(\dfrac{1.8464}{2}\right) - 0.3142(1.8) = 0.6109$.

In general, $\varphi_{rad} = \omega \Delta t = 2\pi f \Delta t$ which can be estimated once the approximate sine function is obtained by the filtering.

———

In summary, the estimation procedure just presented provides a reasonably accurate goodness of fit solution of the real variable sinusoidal function, $x(t) = 2\sin\left(\dfrac{2\pi}{20}t + 35°\right)$, the basis function that generated the noise corrupted time series. We estimate the overall error rate at about 2-3%. The general method will be useful for modeling the parameters of simple harmonic motion degraded by noise.

Calculations for the derived hypothetical sinusoid will show the following structure in tabular and graphical format. See Table 1 below.

Table 1. Structural properties of time-delayed sinusoid, $x(t) = 2\sin\left(\dfrac{2\pi}{20}t + \dfrac{7}{36}\pi\right)$

| Time $t$ | Value at $t$ ($\pi = 3.141593....$) |
|---|---|
| $-1\dfrac{17}{18}$ | $t_0$ (right horizontal shift factor). This value is the time delay, $\Delta t = \dfrac{\varphi}{\omega}$, and represents the lower bound on the graph for one period $\left(T = \dfrac{2\pi}{\omega}\right)$, $\left[\left(-\dfrac{\varphi}{\omega},0\right),\left(-\dfrac{\varphi}{\omega}+\dfrac{2\pi}{\omega},0\right)\right]$, for the function, $y = A\sin(\omega t + \varphi), t \geq \Delta t$. The scaled sine function is created by setting $t \to t \pm \Delta t$ to get $y = A\sin(\omega t)$ with graph points $(0,0),\left(\dfrac{2\pi}{\omega},0\right)$ for one period. Thus, $A\int_{t_a+\Delta t}^{t_b+\Delta t}\sin[\omega(t-\Delta t)+\varphi]dt = A\int_{\frac{a\pi}{\omega}}^{\frac{b\pi}{\omega}}\sin(\omega t)dt$, $t_k = \dfrac{k\pi - \varphi}{\omega}$. $x(t) = 0$ |
| 0 | $1.1472 = A\sin\varphi$ |
| $\pm 2$ | Amplitude $A$ |
| $3\dfrac{1}{18}$ | $t_{\frac{\pi}{2}}$ (max.) $x(t) = +2$ |
| $8\dfrac{1}{18}$ | $t_\pi$ (ist 0–crossover, $t \geq 0$) $x(t) = 0$ |
| $13\dfrac{1}{18}$ | $t_{\frac{3}{2}\pi}$ (min.) $x(t) = -2$ |
| $18\dfrac{1}{18}$ | $t_{2\pi}$ (2nd 0–crossover, $t \geq 0$) $x(t) = 0$ |
| $23\dfrac{1}{18}$ | $t_{\frac{5}{2}\pi}$ (max.); one period on $t \geq 0$ is calculated for time-delayed sinusoids[3] (simple harmonic motion): $T = \dfrac{2\pi}{\omega} = \overbrace{\left(t_\pi - t_{\frac{\pi}{2}}\right)}^{5=\frac{T}{4}} + \overbrace{\left(t_{2\pi} - t_\pi\right)}^{10=\frac{T}{2}} + \overbrace{\left(t_{\frac{5}{2}\pi} - t_{2\pi}\right)}^{5=\frac{T}{4}} = t_{\frac{5}{2}\pi} - t_{\frac{\pi}{2}}$ $= \left(\dfrac{5\pi - 2\varphi}{2\omega}\right) - \left(\dfrac{\pi - 2\varphi}{2\omega}\right) = 23\dfrac{1}{18} - 3\dfrac{1}{18} = 20$. Derived from symmetric odd function definition, 1 period: $\int_a^b \sin(\omega t + \varphi)dt = 0, b - a = \dfrac{2\pi}{\omega}, a \geq 0, \int_{t_{\frac{\pi}{2}}}^{t_{\frac{5\pi}{2}}}\sin(\omega t + \varphi)dt = \int_0^{\frac{2\pi}{\omega}}\sin(\omega t)dt = 0$ and |

---

[3] When phase is negative ($\Delta t > 0$) [time ahead sinusoid, $A\sin(\omega t - \varphi)$], we take $T = t_{2\pi} - t_0 = t_{2\pi} - \Delta t$ with graph points $\left[\left(\dfrac{\varphi}{\omega},0\right),\left(\dfrac{\varphi}{\omega}+\dfrac{2\pi}{\omega},0\right)\right]$.

| | property, $\cos\left(k\dfrac{\pi}{2}\right) = 0, k = 1,3,5,...$ |
|---|---|
| | $x(t) = +2$ |

Example: the value of $t$ in $A\sin(\omega t + \varphi)$ that gives $\pi$ $(t_\pi)$ is found by solving $t$ in $\omega t + \phi = \pi$,

so that $x(t) = 2\sin\left(\dfrac{2\pi}{20}t + 35°\right)\bigg|_{t=8\frac{1}{18}} = 0.$

The following graph shows the hypothetical sinusoid and the transformed sine wave.

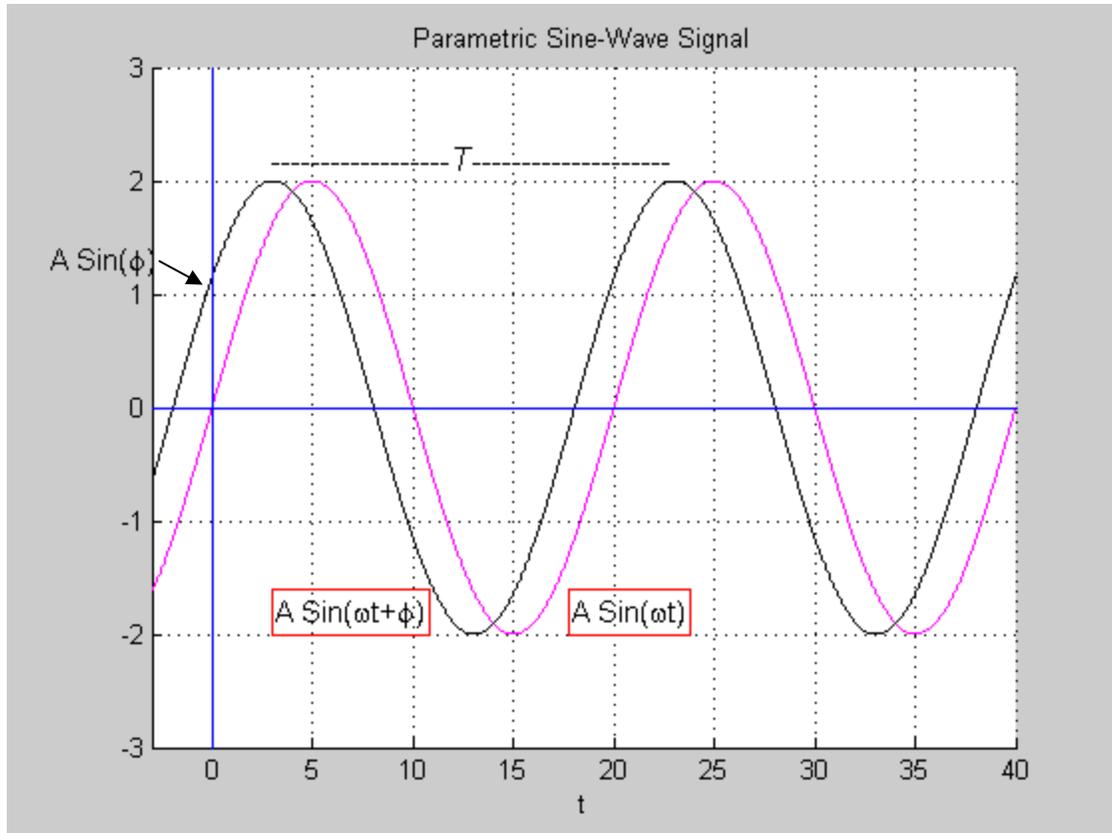

Figure 5. This graph shows three distinct reference points for one period $T = \dfrac{2\pi}{\omega}$:

$\left[\left(-\dfrac{\varphi}{\omega}, 0\right), \left(-\dfrac{\varphi}{\omega} + \dfrac{2\pi}{\omega}, 0\right)\right]$   [periodic function $y = A\sin(\omega t + \varphi), t \geq \dfrac{\varphi}{\omega}$]

$\left[\left(t_{\frac{\pi}{2}}, A\right), \left(t_{\frac{5\pi}{2}}, A\right)\right]$   [simple harmonic motion, $y = A\sin(\omega t + \varphi), t \geq 0$]

$\left[(0,0), \left(\dfrac{2\pi}{\omega}, 0\right)\right]$   [SHM scaled, $y = A\sin(\omega t), t \geq 0$]

The integral of each odd function evaluates to 0.

Note on Derivations

All mathematical solutions were checked on Mathematica. Solutions provided by hand used the following basic identities and elementary relations for the sine and cosine circular functions:

- $\sin A \sin B = \dfrac{1}{2}[\cos(A-B) - \cos(A+B)]$
- $\sin A - \sin B = 2\cos\dfrac{A+B}{2}\sin\dfrac{A-B}{2}$
- $\sin(A-B) = \sin(A)\cos(B) - \sin(B)\cos(A)$
- $\displaystyle\int_u^v \sin(ax+b)\sin(ax+d)\,dx = \dfrac{v-u}{2}\cos(b-d) - \dfrac{\sin[a(v-u)]\cos[a(u+v)+b+d]}{2a}$

## 4. Summary

In conclusion, use of the statistical data smoothing FIR filters and minimization algorithm (8) will provide a goodness of fit empirical curve which allows estimation of all parameters in a sinusoidal model describing simple harmonic motion, $A\sin(\omega t + \varphi)$, for modeling noisy input data $X(t)$ in a finite time series. The method is conducted in two stages. The standard curve fitting approaches are not always successful due to initial value sensitivity and other issues. For these reasons the authors have derived the present method for a range of amplitudes, frequencies, phases and sample sizes that is not dependent on initial value specification.

## 5. Alternatives

The engineering method as described can be used for time–continuous data models as well. Different sinusoidal models such as $x(t) = A\cos(\omega t + \varphi)A\sin[t\omega + \varphi]$ or more complicated sinusoidal models such as $x(t) = A\sin(\omega_1 t + \varphi) + B\cos(\omega_2 t + \varphi)$ are also amenable to the present method based on well known trigonometric reduction techniques (principle of superposition). Also, additional screening tests aside from the Wald–Wolfowitz Runs Test (Wald and Wolfowitz, 1940) and circular correlation can and should be conducted to improve accuracy of the signal–noise decision.

Frequency, the key parameter in the estimation method, is best obtained directly from the Fast Fourier Transform (FFT) of the data set when the noise level is not too high. Frequency can also be estimated from the smoothed data by observing the value corresponding to one period $T$. Then $f = 1/T$; $f$ can also be estimated by observing the one period mark of the circular ACF of ACF of the sinusoid, Eq. (5c).

# APPENDIX I
# Method Steps

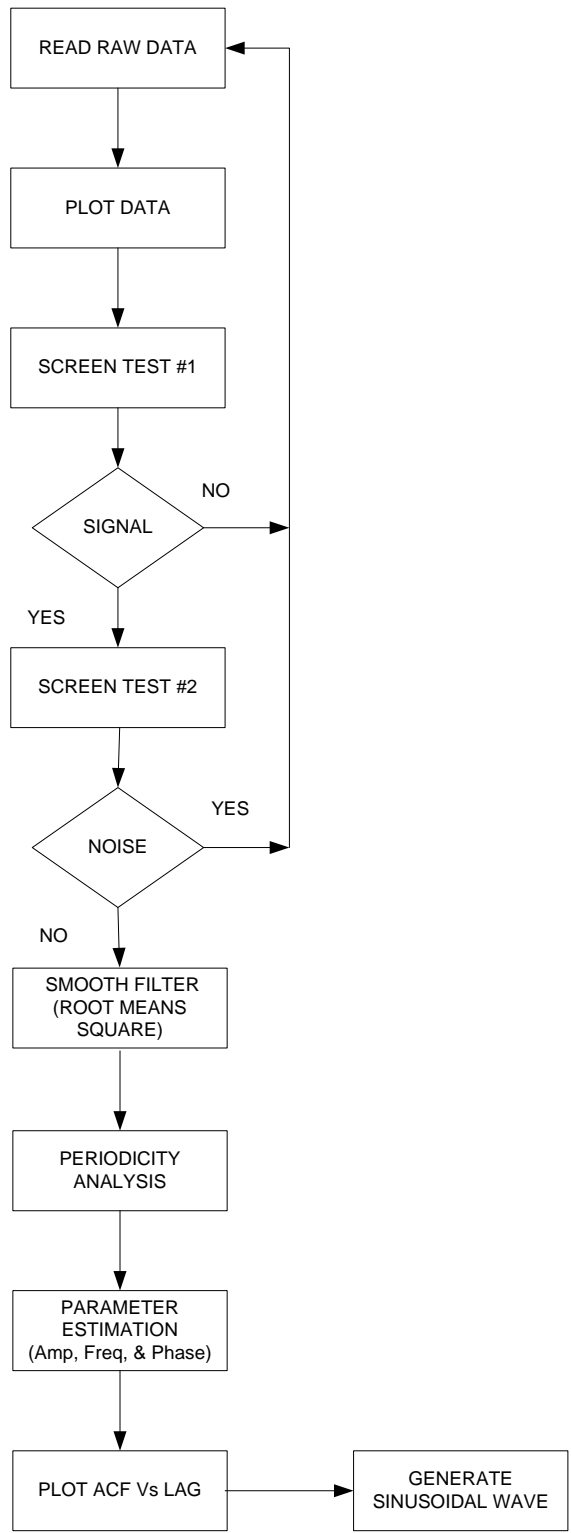

FLOW CHART PROCESS